\documentclass[reqno,12pt]{amsart}
\usepackage{a4wide}

\numberwithin{equation}{section}

\newtheorem{Theorem}{Theorem}
\newtheorem{conj}{Conjecture}

\newtheorem{Proposition}{Proposition}
\newtheorem{Lemma}{Lemma}
\theoremstyle{remark}
\newtheorem*{Remark}{Remark}

\def\al{\alpha}
\def\be{\beta}

\def\ep{\varepsilon}

\def\et{\eta}

\def\Ga{\Gamma}

\def\({\left(}
\def\){\right)}
\def\[{\left[}
\def\]{\right]}
\def\ii{\infty}

\def\cl#1{\left\lceil#1\right\rceil}

\def\dd{\textup{d}}

\def\lcm{\operatorname{lcm}}
\def\Re{\operatorname{\mathfrak R\mathfrak e}}

\begin{document}

\title[]{
An identity of Andrews, multiple integrals,
and very-well-poised hypergeometric series}
\author[]{C. Krattenthaler$^\dagger$ and T. Rivoal}

\address{
Institut Girard Desargues,
Universit{\'e} Claude Ber\-nard Lyon-I,
21, avenue Claude Ber\-nard,
F-69622 Villeurbanne Cedex, France}
\email{kratt@igd.univ-lyon1.fr}
\address{
Laboratoire de Math{\'e}matiques Nicolas Oresme, CNRS UMR 6139,
Universit{\'e} de Caen,  BP 5186,
14032 Caen cedex,
France}
\email{rivoal@math.unicaen.fr}
\thanks{$^\dagger$Research partially supported by 
the programme  ``Improving the Human Research Potential'' 
of the European Commission, grant HPRN-CT-2001-00272, 
``Algebraic Combinatorics in Europe"}



\begin{abstract}
We give a new proof of a theorem of Zudilin that equates a
very-well-poised hypergeometric series and a particular multiple integral.
This integral generalizes integrals of Vasilenko and Vasilyev 
which were proposed as tools in the study of the arithmetic behaviour of 
values of the Riemann zeta function at integers.
Our proof is based on limiting cases of a basic 
hypergeometric identity of Andrews. 
\end{abstract}

\maketitle


\section{Introduction}
After Ap{\'e}ry's 1978 proof of the irrationality of $\zeta(2)$ and
$\zeta(3)$ (see \cite{ap}), $\zeta(s)$ denoting the Riemann zeta
function, Beukers
\cite{be}  
gave another proof  with the help of his famous integrals
$$
\int_{[0,1]^2}\frac{x^n(1-x)^ny^n(1-y)^n}{(1-(1-x)y)^{n+1}}
\,\dd x\,\dd y =\al_n\zeta(2)-\be_n
$$
and
$$
\int_{[0,1]^3}\frac{x^n(1-x)^ny^n(1-y)^nz^n(1-z)^n}
{(1-(1-(1-x)y)z)^{n+1}}\,\dd x\,\dd y\,\dd z
=2a_n\zeta(3)-b_n.
$$
Here, $n\ge 0$ is an integer, and $\al_n$, $\be_n$, $a_n$, $b_n$
are rational numbers. More precisely, $\al_n$, $a_n$, 
$\textup{d}_n^2 \be_n$ and $\textup{d}_n^3 b_n$ are integers, with 
$\textup{d}_n=\lcm\{1,2, \ldots, n\}$. 

Extending Vasilenko's method \cite{vas}, Vasilyev \cite{va1, va2}  
considered a family of integrals  
generalizing Beukers' pattern:
\begin{equation} \label{eq:vas} 
J_{E,n}=\int_{[0,1]^E}\frac{\prod_{i=1}^E x_i^n(1-x_i)^n }
{Q_E(x_1,x_2,\ldots, x_E)^{n+1}}\,\dd x_1\,\dd x_2\cdots \dd x_E,
\end{equation}
where 
\begin{equation*}
Q_E(x_1,x_2,\ldots, x_E)=
1-(\cdots (1-(1-x_E)x_{E-1})\cdots)x_1,
\end{equation*}
and where $E$ is an integer $\ge 2$. 

He then formulated the following conjecture, which he proved 
for $E=4, 5$  and which is also 
true for $E=2, 3$ because of Beukers' work.
\begin{conj}[Vasilyev]
\label{conj3}\leavevmode\newline
{\em (i)} For all integers $E\ge 2$ and $n\ge 0$, there exist 
rational numbers $p_{m,E,n}$ such that
\begin{equation}\label{eq:vasint}
J_{E,n}=p_{0,E,n}+\underset{m \equiv E \,
(\textup{mod}\, 2)}{\sum_{m=2, \ldots, E}} p_{m,E,n} \zeta(m).
\end{equation}
{\em (ii)} Furthermore, $\textup{d}_n^{E}p_{m,E,n}$ is an integer for all 
$m=0, 1, \ldots, E$.
\end{conj}
Part (i) of this conjecture has been proved by Zudilin in \cite[Sec.~8]{zu2}, 
thanks to an unexpected identity between a certain multiple integral $J_m$ 
(generalizing Vasilyev's ones) and a non-terminating
very-well-poised hypergeometric series (see Theorem~\ref{thm:zudgene} below). 
Part (ii) was also proved by Zudilin in \cite{zu2}, except
for the coefficient $p_{0,E,n}$. 
A  sharper version of Part~(ii)
has been established by the authors \cite{kratriv} 
when $E$ is odd (including the claim about $p_{0,E,n}$), 
and for all coefficients but $p_{0,E,n}$ when $E$ is even. 

The aim of this note is to give a new proof of Zudilin's identity, which is the
content of Theorem~\ref{thm:zudgene} in the next section.  
In fact, our proof shows that, modulo a more or less evident 
expansion of the Vasilyev-type integral $J_m$
as a multiple sum
(see Proposition~\ref{thm:zlobin} in Section~\ref{sec:int=andr}), 
Zudilin's identity is a limiting
case of a thirty year old identity between a terminating multiple sum and 
a terminating very-well-poised hypergeometric series due to Andrews 
(see \eqref{eq:andrews} below).
We believe that this is an interesting observation because attempts to
prove Zudilin's identity by manipulating hypergeometric series
directly failed because
of convergence problems. Zudilin circumvents these problems by having recourse
to a Barnes-type (i.e., contour) integral in place of the
very-well-poised hypergeometric series.
Our proof shows that there is indeed
a ``purely hypergeometric" proof (i.e., a proof just using summation and
transformation formulas for hypergeometric series), 
but to be able to accomplish it, one has
to go ``one level higher in hierarchy," meaning that one finds a 
{\it terminating}
identity ``above," of which the identity which one actually wants to prove
is a limiting case.\footnote{This point is also of interest ``philosophically."
There are several proposers (of whom Koornwinder \cite{Koorn} 
seems to have been
the first; see \cite[Remark~3.2]{Schloss} and \cite[paragraph after
Eq.~({\it Apery})]{ZeilAQ} for printed versions) 
of the ``conjecture" that above
every identity for non-terminating hypergeometric series (which are very often
difficult to prove; in particular, the automatic tools described in 
\cite{PeWZAA} 
do not apply) there sits a more general identity for terminating
series
(which, at least in principle, can be proved automatically), of which the
non-terminating identity is a limiting case. Of course, the analyst would 
object
that above every terminating identity there exists an even more general
non-terminating identity, of which the terminating one is a special case. 
Clearly, this dispute is as easy to settle
as the dispute about the question of which was first, the hen or the egg \dots} 
This identity ``above" is Andrews' identity, and it does
indeed have a purely hypergeometric proof (see Section~\ref{sec:andrews}
for more information).

As an aside, we mention that, in the recent paper \cite{zlo}, Zlobin
shows that the multiple integral $J_m$ is also equal to an
integral of the same type as those of Sorokin in \cite{so,so1} 
(see also Fischler~\cite{fiscvar} for similar results). 
The latter integral can also be expanded as a multisum, in a
manner completely analogous to the way we derive the multisum
expansion for $J_m$ in the proof of
Proposition~\ref{thm:zlobin}. As a result, one obtains again exactly the
right-hand sides of \eqref{eq:Zlobin1} and \eqref{eq:Zlobin2},
respectively. Thus, this provides an alternative proof of Zlobin's
result. As a matter of fact, when this work originated, we went the other way,
that is, our starting point was the multisum expansion of Zlobin's
integral, until we realized that, actually, the integral $J_m$ admits
the same treatment.

Zudilin's identity is recalled in Theorem~\ref{thm:zudgene}
in the next section. 
The limiting cases of Andrews' identity which we need are stated and proved
in Proposition~\ref{prop:even} in Section~\ref{sec:andrews}. 
One of the lemmas which we need for carrying out these limits
generalizes a lemma of Zhao \cite{ZhaoAB} on the convergence of
multizeta functions, see the remark after the proof of Lemma~\ref{lem:4}.
The purpose of Section~\ref{sec:int=andr} is to relate these
identities to the Vasilyev-type integral $J_m$, 
see Proposition~\ref{thm:zlobin}.
We finally prove Theorem~\ref{thm:zudgene}  in
Section~\ref{sec:proofvaszud}. 

\section{Zudilin's identity}\label{sec:zud}

In order to be able to state Zudilin's identity,
we need to recall the standard notation for (generalized) 
hypergeometric series,
\begin{eqnarray*}
{}_{p+1}F_p\left[\begin{matrix}
\alpha_0,\alpha_1,\ldots,\alpha_{p}\\ 
\beta_1,\ldots,\beta_p\end{matrix};z\right]
=\sum_{k=0}^{\ii}
\frac{(\alpha_0)_k\,(\alpha_1)_k\cdots(\alpha_{p})_k}
{k!\,(\beta_1)_k\cdots(\beta_p)_k} \,z^k,
\label{eq:hyper}
\end{eqnarray*}
where $p\ge 1$, 
$\alpha_j\in\mathbb{C}$, 
$\beta_j\in\mathbb{C}\setminus\mathbb{Z}_{\le 0}$ and, by definition,  
$(x)_0=1$ and $(x)_{\ell}=x(x+1)\cdots(x+\ell-1)$ for $\ell\ge1$.
The series is absolutely convergent for all $z\in\mathbb{C}$ 
such that $\vert z\vert <1$,
and also for $ \vert z \vert =1$ provided 
$\Re(\be_1+\dots+\be_p)>\Re(\al_0+\al_1+\dots+\al_p)$.
Furthermore, it is said to be {\it balanced} if 
$\al_0+\cdots+\al_p+1=\be_1+\cdots +\beta_p$~and {\it
  very-well-poised} if $\al_0+1=\al_1+\be_1=\cdots=\al_p+\be_p$ and
$\alpha_1=\frac12\,\alpha_0+1$. See the books \cite{AAR, bail,
  GaRaAA, SlatAC} for more information on hypergeometric series.

Let $z$, $a_0, a_1,\dots,a_{m}$, and $b_1,\dots,b_m$ be complex numbers
such that $\vert z\vert <1$, $\Re(b_i)>\Re(a_i)>0$ 
for all $i=1,2,\dots,m$,
and let us define the Vasilyev-type integral
\begin{equation}\label{eq:Jdef}
J_{m}\left[\begin{matrix}
a_0,a_1,\dots,a_m\\b_1,\dots,b_m\end{matrix};z\right]
=\int _{[0,1]^m} ^{}\frac {\prod _{i=1}
^{m}x_i^{a_i-1}(1-x_i)^{b_i-a_i-1}}
{(1-(1-(\cdots(1-x_m)x_{m-1})\cdots)x_1z)^{a_0}}
\,\dd x_1\,\dd x_2\cdots \dd x_m,
\end{equation}
which is absolutely convergent under the above conditions. 
(We will sometimes
use the short notation $J_m$ for this integral if there is no ambiguity
about the parameters.) 
It is also absolutely convergent for $z=1$, 
provided that we also assume that $\Re(b_1-a_1)>\Re(a_0)$ if $m=1$,
respectively  $\Re(b_1-a_1)\ge\Re(a_0)$ if $m>1$. 
Since previous authors assume more restrictive conditions 
in the case $z=1$ (in
particular, restrictions that are not satisfied by Vasilyev's
integrals \eqref{eq:vas}), 
we sketch the verification of the convergence here for the sake
of completeness. 

If
$m=1$, then $J_m=J_1$ is a beta integral. If $m\ge2$, then, because of
$$ 1-(1-(\cdots(1-x_m)x_{m-1})\cdots)x_1\ge 1-x_1$$ 
and 
$\Re(b_1-a_1)\ge\Re(a_0)$, we have
\begin{align*}
\int _{\ep_m} ^{1-\ep_m}&
\dots\int _{\ep_2} ^{1-\ep_2}
\int _{\ep_1} ^{1-\ep_1}\left\vert\frac {\prod _{i=1}
^{m}x_i^{a_i-1}(1-x_i)^{b_i-a_i-1}}
{(1-(1-(\cdots(1-x_m)x_{m-1})\cdots) x_1)^{a_0}}
\right\vert
\,\dd x_1\,\dd x_2\cdots \dd x_m\\
&\le
\int _{\ep_m} ^{1-\ep_m}
\dots\int _{\ep_2} ^{1-\ep_2}
\int _{\ep_1} ^{1-\ep_1}\frac {\prod _{i=2}
^{m}x_i^{\Re(a_i-1)}(1-x_i)^{\Re(b_i-a_i-1)}}
{(1-x_1+Xx_2 x_1)}
\,\dd x_1\,\dd x_2\cdots \dd x_m\\
&\le
\int _{\ep_m} ^{1-\ep_m}\dots
\int _{\ep_2} ^{1-\ep_2} \left({\prod _{i=2}
^{m}x_i^{\Re(a_i-1)}(1-x_i)^{\Re(b_i-a_i-1)}}\right)\\
&\kern4cm
\cdot
\left(-\frac {\log(1-x_1+Xx_2x_1)} 
{1-Xx_2}\right)\bigg\vert_{x_1=\ep_1}^{1-\ep_1}
\,\dd x_2\cdots \dd x_m,
\end{align*}
for any small $\ep_1,\ep_2,\dots,\ep_m>0$, where we wrote $X$ for
$1-(\cdots(1-x_m)x_{m-1})\cdots)x_3$. (In case that $m=2$, $X$ has to
be interpreted as 1.) If we perform the limit
$\ep_1\to0$, then the right-hand side of this inequality becomes the integral
\begin{equation} \label{eq:int2} 
\int _{\ep_m} ^{1-\ep_m}\dots
\int _{\ep_2} ^{1-\ep_2} \left({\prod _{i=2}
^{m}x_i^{\Re(a_i-1)}(1-x_i)^{\Re(b_i-a_i-1)}}\right)\\
\left(-\frac {\log(Xx_2)} 
{1-Xx_2}\right)
\,\dd x_2\cdots \dd x_m.
\end{equation}
In the integrand, there is no problem as $Xx_2\to1$, since the
function
$\log(Xx_2)/(1-Xx_2)$ is continuous at $Xx_2=1$. On the other hand,
if we fix $\et>0$, then
for $Xx_2$ sufficiently close to 0, we have 
$$\vert\log(Xx_2)\vert<(Xx_2)^{-\et}\le (1-x_3)^{-\et}x_2^{-\et}.$$
Thus, choosing $\et=\frac {1} {2}\min\{\Re(a_2),\Re(b_3-a_3)\}$, we
see that the integral in \eqref{eq:int2}, and thus the original
integral $J_m$, exists.

\begin{Theorem}[Zudilin]\label{thm:zudgene} For every integer $m\ge 1$,
  the following identity holds:
\begin{multline} \label{eq:Zudilin}
J_{m}\left[\begin{matrix}
h_1,h_2, h_3, \ldots,h_{m+1}\\ 
1+h_0-h_3, 1+h_0-h_4, \ldots, 1+h_0-h_{m+2}\end{matrix};1\right]
\\
=
\frac{\Ga(1+h_0)\prod_{j=3}^{m+1}\Ga(h_j)}{\prod_{j=1}^{m+2}\Ga(1+h_0-h_j)}\cdot
\left(\prod_{j=1}^{m+1}\Ga(1+h_0-h_j-h_{j+1})\right)\kern5cm\\
\times
\;{}_{m+4}F_{m+3}
\left[\begin{matrix}h_0, \frac{1}{2}h_0+1, h_1, \ldots,  h_{m+2}\\ 
\frac{1}{2}h_0, 1+h_0-h_1, \ldots,  1+h_0-h_{m+2}
\end{matrix};(-1)^{m+1}\right],
\end{multline}
provided that 
$1+\Re(h_0)>\frac{2}{m+1}\sum_{j=1}^{m+2}\Re(h_j)$,
$\Re(1+h_0-h_{j+1})>\Re(h_j)>0$ for $j=2,3,\dots,m+1$,
and 
$\Re(1+h_0-h_{3}-h_2)\ge\Re(h_1)$, 
these conditions 
ensuring that both sides of~\eqref{eq:Zudilin} are well-defined.
\end{Theorem}

In the case of the original 
integrals $J_{E,n}$  of Vasilyev, the identity in
Theorem~\ref{thm:zudgene} reads as follows: 
for any integers $n\ge 0$ and $E\ge 2$,  
\begin{equation*}\label{vaszud}
J_{E,n}=\frac{n!^{2E+1}(3n+2)!}{(2n+1)!^{E+2}}
\;{}_{E+4}F_{E+3}
\left[\begin{matrix}3n+2, \frac{3}{2}n+2, n+1, \ldots,   n+1 \\ 
\frac{3}{2}n+1, 2n+2, \ldots, 2n+2
\end{matrix};(-1)^{E+1}\right].
\end{equation*} 
{}From \cite{br, ri1}, it follows that such a very-well-poised 
hypergeometric series gives rise to a decomposition of the shape
\eqref{eq:vasint}.

\section{Limiting cases of Andrews' hypergeometric identity}\label{sec:andrews}
Let $N$ and $s$ be positive integers, and $a$, $b_1, \ldots, b_{s+1}$,
$c_1, \ldots, c_{s+1}$ 
be complex numbers such 
that none of $1+a-b_j$,  $1+a-c_j$, $j=1,2, \ldots, s+1$, and 
$1+a+N$ are non-positive integers.

Andrews' identity \cite[Theorem~4]{Andr} relates a terminating 
very-well-poised basic hypergeometric series
to a terminating multiple basic hypergeometric series. 
We shall need here the limiting case of this identity 
when $q\to1$, so that the series
there reduce to ``ordinary" hypergeometric series.
That is, we replace $a$ by $q^a$, $b_i$ by $q^{b_i}$, $c_i$ by $q^{c_i}$, 
there, 
and then let $q$ tend to 1. The result can be compactly written in the form
\begin{multline}
{}_{2s+5}F_{2s+4}\left[
\begin{array}{c}
a,\frac{a}{2}+1, b_1,  c_1, \ldots,  b_{s+1},  c_{s+1}, -N \\
     \frac{a}{2}, 1+a-b_1, 1+a-c_1, \ldots,  1+a-b_{s+1}, 1+a-c_{s+1}, 1+a+N 
\end{array}
\,;\, 1 \right]
\\
=
\frac{(1+a)_N\,(1+a-b_{s+1}-c_{s+1})_N}{(1+a-b_{s+1})_N\,(1+a-c_{s+1})_N}
\sum_{k_1, k_2,\ldots, k_{s} \ge 0} 
\frac{(-N)_{k_1+\cdots
+k_s}}{(b_{s+1}+c_{s+1}-a-N)_{k_1+\cdots+k_s}}\kern2cm\\
\cdot
\prod_{j=1}^s
\frac{(1+a-b_{j}-c_j)_{k_j}\,(b_{j+1})_{k_1+\cdots+k_j}\,(c_{j+1})_{k_1+\cdots+k_j}}
{k_j!\,(1+a-b_j)_{k_1+\cdots+k_j}\,(1+a-c_j)_{k_1+\cdots+k_j}}.
\label{eq:andrews}
\end{multline}
The proof in \cite{Andr} uses
Whipple's transformation between a balanced $_4F_3$-series and a
very-well-poised $_7F_6$-series, 
\begin{multline*}
{} _{4} F _{3} \!\left [ \begin{matrix} { a, b, c, -N}\\
{ e, f, 1 + a + b + c - e -
   f - N}\end{matrix} ; {\displaystyle 1}\right ] =
\frac {( -a - b + e + f)_N\,( -a - c + e
+ f) _{N}} {( -a + e + f)_N\, (-a - b - c + e + f) _{N}} \\
\cdot
  {} _{7} F _{6} \!\left [ \begin{matrix} { -1 - a + e + f, \frac{1}{2} -
\frac{a}{2} + \frac{e}{2} + \frac{f}{2}, -a + f, -a + e, b, c, -N}\\ {
- \frac{1}{2}  - \frac{a}{2} + \frac{e}{2} + \frac{f}{2}, e,
f, -a - b + e + f, -a - c + e + f, -a + e + f + N}\end{matrix}\; ;
{\displaystyle 1}\right ],
\end{multline*}
and the Pfaff--Saalsch{\"u}tz summation in an iterative fashion.
In particular, the identity \eqref{eq:andrews} reduces to Whipple's
transformation for $s=1$.

We prove that the same kind of identity holds for {\it non-terminating} 
hypergeometric series provided the parameters $a$, $b_j$, and $c_j$,
$j=1,2,\dots,s+1$, 
satisfy some further conditions.

\begin{Proposition} \label{prop:even}
{\em(i)} Let $s\ge 1$ be an integer, and let $a$, $b_1, \ldots, b_{s+1}$, $c_1,
\ldots, c_{s+1}$ 
be complex numbers such 
that none of\/ $1+a-b_j$,  $1+a-c_j$, $j=1,2, \ldots, s+1$, is a
non-positive integer. Furthermore, we assume that 
\begin{equation} \label{eq:cond1} 
\Re\left((2s+1)(a+1)-
2\sum _{j=1} ^{s+1}(b_j+c_j)\right)>0
\end{equation}
and
\begin{equation} \label{eq:cond2}
\Re\left((1+a-b_{s+1}-c_{s+1})+\sum _{j=r} ^{s}A_j(1+a-b_j-c_j) \right)>0
\end{equation}
for all $r=2,3,\dots,s+1$ {\em(}in the case that $r=s+1$, the empty sum
$\sum _{j=r} ^{s}$ has to be interpreted as $0${\em)}, for all
possible choices 
of $A_j=1$ or $2$, for $j=2,3,\dots,s$. Then 
\begin{multline}
{}_{2s+4}F_{2s+3}\left[
\begin{array}{c}
a,\frac{a}{2}+1, b_1,  c_1, \ldots,  b_{s+1},  c_{s+1} \\
     \frac{a}{2}, 1+a-b_1, 1+a-c_1, \ldots,  1+a-b_{s+1}, 1+a-c_{s+1}
\end{array}
\,;\, -1 \right]
\\
=
\frac{\Ga(1+a-b_{s+1})\,\Ga(1+a-c_{s+1})}{\Ga(1+a)\,\Ga(1+a-b_{s+1}-c_{s+1})}
\kern6cm\\
\times
\sum_{k_1, k_2,\ldots, k_{s} \ge 0}
\prod_{j=1}^s
\frac{(1+a-b_{j}-c_j)_{k_j}\,(b_{j+1})_{k_1+\cdots+k_j}\,(c_{j+1})_{k_1+\cdots+k_j}}
{k_j!\,(1+a-b_j)_{k_1+\cdots+k_j}\,(1+a-c_j)_{k_1+\cdots+k_j}}.
\label{eq:Feven}
\end{multline}

{\em(ii)} Let $s\ge 1$ be an integer, and let $a$, $c_0$, 
$b_1, \ldots, b_{s}$, $c_1, \ldots, c_{s}$ be complex numbers such 
that none of\/
$1+a-b_j$,  $1+a-c_j$, $j=0,1, \ldots, s$, is a non-positive
integer. Furthermore, we assume that 
\begin{equation} \label{eq:cond1a} 
\Re\left(2s(a+1)-2c_0-
2\sum _{j=1} ^{s}(b_j+c_j)\right)>0,
\end{equation}
\begin{equation} \label{eq:cond2a}
\Re\left((1+a-b_{s}-c_{s})+\sum _{j=r} ^{s-1}A_j(1+a-b_j-c_j) \right)>0
\end{equation}
for all $r=2,3,\dots,s$ {\em(}in the case that $r=s$, the empty sum
$\sum _{j=r} ^{s-1}$ has to be interpreted as $0${\em)}, and
\begin{equation} \label{eq:cond3a}
\Re\left((1+a-c_0-b_{1}-c_{1})+\sum _{j=2} ^{s-1}A_j(1+a-b_j-c_j) \right)>0,
\end{equation} 
for all possible choices of $A_j=1$ or $2$, for $j=2,3,\dots,s-1$. Then 
\begin{multline}
{}_{2s+3}F_{2s+2}\left[
\begin{array}{c}
a,\frac{a}{2}+1, c_0, b_1,  c_1, \ldots,  b_{s},  c_{s} \\
     \frac{a}{2},1+a-c_0,  1+a-b_1, 1+a-b_1, \ldots,  1+a-b_{s}, 1+a-c_{s}
\end{array}
\,;\, 1 \right]
\\
=
\frac{\Ga(1+a-b_{s})\,\Ga(1+a-c_{s})}{\Ga(1+a)\,\Ga(1+a-b_{s}-c_{s})}
\sum_{k_1, k_2,\ldots, k_{s} \ge 0}
\frac{(b_1)_{k_1}\,(c_1)_{k_1}}
{k_1!\,(1+a-c_0)_{k_1}}\kern4cm\\
\cdot
\prod_{j=2}^{s}
\frac{(1+a-b_{j-1}-c_{j-1})_{k_j}\,(b_{j})_{k_1+\cdots+k_j}\,(c_{j})_{k_1+\cdots+k_j}}
{k_j!\,(1+a-b_{j-1})_{k_1+\cdots+k_j}\,(1+a-c_{j-1})_{k_1+\cdots+k_j}}.
\label{eq:Fodd}
\end{multline}
\end{Proposition}

Our proof of this proposition is based on three lemmas, which we state
and prove first.

\begin{Lemma} \label{lem:2}
Let $\alpha$ and $\beta$, $i=1,2,\dots,m$, be complex numbers such that
$\beta$ is not a non-positive integer. Then for any non-negative integer
$k$ we have
$$
\left\vert\frac {\Ga(\alpha+k)} {\Ga(\beta+k)}\right\vert
\le D_1 \cdot(k+1)^{\Re(\alpha-\beta)},$$
where $D_1$ is a constant which does not depend on $k$.
\end{Lemma}
\begin{proof}By Stirling's formula, we have
$$\frac {\Ga(\alpha+k)} {\Ga(\beta+k)}\sim 
(k+1)^{\alpha-\beta}$$
as $k\to\infty$. Hence, the claim follows immediately.
\end{proof}

\begin{Lemma} \label{lem:3}
Let $A$ and $B$ be real numbers such that $A+B+1<0$, and let $C$ be a
non-negative integer. 
Then, for any non-negative integer $h$, we have
$$
\sum _{k=0} ^{\infty}(k+1)^A(h+k+1)^B\left(\log(h+k+2)\right)^C
\le D_2 (h+1)^{\max\{B,A+B+1\}}\left(\log(h+2)\right)^{C+1},$$
where $D_2$ is a constant independent of $h$.
\end{Lemma} 
\begin{proof}
We split the summation range into the ranges
$R_0=\{0,1,\dots,2^{\cl{\log_2(h+1)}+1} -h-1\}$ and
$$R_s=\{2^{s}-h,2^s-h+1,\dots,2^{s+1}-h-1\},\quad \quad s=\cl{\log_2
(h+1)}+1,\cl{\log_2(h+1)}+2,\dots,$$
where $\cl{x}$ denotes the least integer $\ge x$.
Since, depending on whether $A$ and $B$ are positive or not, 
for $k\in R_s$, $s>0$, we have
\begin{equation*}
(k+1)^A(h+k+1)^B\left(\log(h+k+2)\right)^C
\le \left(\log2^{s+2}\right)^C\cdot
\begin{cases} 2^{(s+1)(A+B)}&\text{if }A,B\ge 0,\\
2^{(s-1)A+(s+1)B}&\text{if }A<0,\ B\ge 0,\\
2^{(s+1)A+sB}&\text{if }A\ge 0,\ B<0,\\
2^{(s-1)A+sB}&\text{if }A,B<0,\\
 \end{cases}
\end{equation*}
which implies that
$$(k+1)^A(h+k+1)^B\left(\log(h+k+2)\right)^C
\le D_3 \cdot (s+2)^C\cdot 2^{s(A+B)},$$
with a constant $D_3$ which is independent of $h$. 
Thus, for the sum
over the range $\{k\ge 2^{\cl{\log_2(h+1)}+1} -h\}$ we have
\begin{align*}
\sum _{k= 2^{\cl{\log_2(h+1)}+1} -h} ^{\infty}(k+1)^A&(h+k+1)^B
\left(\log(h+k+2)\right)^C\\
&=
\sum _{s=\cl{\log_2
(h+1)}+1} ^{\infty}
\sum _{k\in R_s} ^{}(k+1)^A(h+k+1)^B\left(\log(h+k+2)\right)^C\\
&\le \sum _{s=\cl{\log_2
(h+1)}+1} ^{\infty}
D_3 \cdot 2^s\cdot(s+1)_{C+1}\cdot 2^{s(A+B)}\\
&\le D_3\cdot2^{(\cl{\log_2
(h+1)}+1)(A+B+1)}\frac {(C+1)!} {\left(1-2^{A+B+1}\right)^{C+2}}\\
&\le D_4\cdot (h+1)^{A+B+1},
\end{align*}
for a constant $D_4$ independent of $h$.

Now we consider the remaining range,
$R_0=\{0,1,\dots,2^{\cl{\log_2(h+1)}+1} -h-1\}$. 
For any $k\in R_0$ we have $k\le 3h+3$, and therefore
$$(h+k+1)^B\left(\log(h+k+2)\right)^C\le \left(\log(4h+5)\right)^C
\cdot
\begin{cases} (h+1)^B&\text{if }B\le0,\\
(4h+4)^B&\text{if }B>0.\end{cases}$$
In particular, there is a constant $D_5$ independent of $k$ such that
$$(h+k+1)^B\left(\log(h+k+2)\right)^C\le D_5
\left(\log(h+2)\right)^C (h+1)^B$$
for all $k\in R_0$. Using this fact, we are able to conclude that
$$
\sum _{k\in R_0} ^{}(k+1)^A(h+k+1)^B\left(\log(h+k+2)\right)^C
\le D_5\cdot\left(\log(h+2)\right)^C
(h+1)^B\sum _{k\in R_0} ^{}(k+1)^A.
$$
Now, if $A<-1$, then $\sum _{k\in R_0} ^{}(k+1)^A<\zeta(-A)$.
If $A=-1$, then $\sum _{k\in R_0} ^{}(k+1)^{-1}<\log(4h+4)$. 
Finally, if $A>-1$, then
$$\sum _{k\in R_0} ^{}(k+1)^A<\int _{0} ^{4h+4}x^A\,\dd x=\frac {1}
{A+1} (4h+4)^{A+1}.$$
In all cases, we obtain that
$$\sum _{k\in R_0} ^{}(k+1)^A(h+k+1)^B\left(\log(h+k+2)\right)^C<
D_6\left(\log(h+2)\right)^{C+1}(h+1)^{\max\{B,A+B+1\}},$$
where $D_6$ is a constant independent of $h$.

To conclude the proof of the lemma, 
the two estimates for the two ranges are combined,
and the
claimed result follows.
\end{proof}

In the statement of the next lemma, we use the following notation:
given two sets $S$ and $T$,
we write $S+T$ for the sum-set $\{x+y:x\in S\text{ and }y\in T\}$. 

\begin{Lemma} \label{lem:4}
Let $E_j$ and $F_j$ be real numbers and let $Z_j$
denote the set $\{F_j,E_j+F_j+1\}$, $j=1,2,\dots,s$.
If 
\begin{equation} \label{eq:bed} 
E_r+F_r+1+\max(Z_{r+1}+Z_{r+2}+\dots+Z_s)<0
\end{equation}
for $r=1,2,\dots,s$, then
the multiple series
\begin{equation} \label{eq:msum} 
\sum _{k_1,\dots,k_s\ge0} ^{}\prod _{j=1} ^{s}(k_j+1)^{E_j}\,
(k_1+\dots+k_j+1)^{F_j}
\end{equation}
converges.
\end{Lemma} 

\begin{proof}
By
applying Lemma~\ref{lem:3} iteratively, we have
\begin{align*}
\sum _{k_1,\dots,k_{s}\ge0} ^{}&\prod _{j=1} ^{s}(k_j+1)^{E_j}\,
(k_1+\dots+k_j+1)^{F_j}\\
&\le
D_7\sum _{k_1,\dots,k_{s-1}\ge0} ^{}\left(\prod _{j=1} ^{s-1}(k_j+1)^{E_j}\,
(k_1+\dots+k_{j}+1)^{F_{j}}\right)\\
&\kern3cm
\cdot(k_1+\dots+k_{s-1}+1)^{\max (Z_s)}\log(k_1+\dots+k_{s-1}+2)\\
&\le
D_8\sum _{k_1,\dots,k_{s-2}\ge0} ^{}\left(\prod _{j=1} ^{s-2}(k_j+1)^{E_j}\,
(k_1+\dots+k_{j}+1)^{F_{j}}\right)\\
&\kern3cm
\cdot(k_1+\dots+k_{s-2}+1)^{\max (Z_{s-1}+Z_s)}
\left(\log(k_1+\dots+k_{s-2}+2)\right)^2,
\end{align*}
and, after the $t$-th iteration, $1\le t\le s-1$,
\begin{align*}
\sum _{k_1,\dots,k_{s}\ge0} ^{}&\prod _{j=1} ^{s}(k_j+1)^{E_j}\,
(k_1+\dots+k_j+1)^{F_j}\\
&\le
D_9\sum _{k_1,\dots,k_{s-t}\ge0} ^{}\left(\prod _{j=1} ^{s-t}(k_j+1)^{E_j}\,
(k_1+\dots+k_{j}+1)^{F_{j}}\right)\\
&\kern2.5cm
\cdot
(k_1+\dots+k_{s-t}+1)^{\max (Z_{s-t+1}+\dots+Z_s)}
\left(\log(k_1+\dots+k_{s-t}+2)\right)^{t}.
\end{align*}
To justify these steps, we have to verify that the condition $A+B+1<0$ 
in Lemma~\ref{lem:3} is satisfied in each iteration. However, this is
exactly the condition \eqref{eq:bed} with $r$ replaced by $s-t$.

Thus, for $t=s-1$ we arrive at the estimate
\begin{align*}
\sum _{k_1,\dots,k_{s}\ge0} ^{}\prod _{j=1} ^{s}(k_j+1)^{E_j}&\,
(k_1+\dots+k_j+1)^{F_j}\\
&\le
D_{10}\sum _{k_1\ge0} ^{}(k_1+1)^{E_1+F_1+\max (Z_{2}+\dots+Z_s)}
\left(\log(k_1+2)\right)^{s-1}.
\end{align*}
Since the sum over $k_1$ at the right-hand side 
converges because of \eqref{eq:bed} with
$r=1$, the claim follows.
\end{proof}

\begin{Remark}
A careful check of our arguments reveals that, in fact, the conditions
in Lemma~\ref{lem:4} are optimal, meaning that they describe {\it
exactly} the domain of convergence of the multiple sum
\eqref{eq:msum}. This can be seen by verifying that, if condition
\eqref{eq:bed} is violated for a particular $r$, then the subsum 
\begin{equation*} 
\sum _{k_r,\dots,k_s\ge0} ^{}\prod _{j=1} ^{s}(k_j+1)^{E_j}\,
(k_1+\dots+k_j+1)^{F_j}
\end{equation*}
of \eqref{eq:msum} does not converge.
Thus, this lemma generalizes Proposition~1
in \cite{ZhaoAB}. It does at the same time {\it correct} that proposition,
and it answers the question raised after the (incomplete) proof of the
proposition. The question, which is asked there, is to determine the
domain of absolute convergence of the multizeta function
\begin{equation} \label{eq:multizeta} 
\zeta(s_d,s_{d-1},\dots,s_1)=
\sum _{0<n_1<\dots<n_d} ^{}\frac {1} {n_1^{s_1}n_2^{s_2}\cdots
n_d^{s_d}}.
\end{equation}
Proposition~1 in \cite{ZhaoAB} states that,
for all $d$-tuples $(s_1,s_2,\dots,s_d)$ with $\Re(s_d)>1$ and 
$\sum _{i=1} ^{d}\Re(s_i)>d$, the series $\zeta(s_d,s_{d-1},\dots,s_1)$
converges absolutely. (As the case $d=3$, $s_1=3$, $s_2=-1$, $s_3=2$ 
shows, these conditions are not sufficient.)

Applying Lemma~\ref{lem:4} to
$$\sum _{0<n_1\le \dots\le n_d} ^{}\frac {1} {n_1^{s_1}n_2^{s_2}\cdots
n_d^{s_d}}=
\sum _{k_1,\dots,k_d\ge0} ^{}
\prod_{j=1}^d \frac 1 {(k_1+k_2+\dots+k_j+1)^{s_j}}$$
(that is, one chooses $s=d$, $E_i=0$ and 
$F_i=-\Re(s_i)$, $i=1,2,\dots,d$, there),
it is seen that the domain of absolute convergence of 
this latter multisum is the
set of all $d$-tuples $(s_1,s_2,\dots,s_d)$ such that 
\begin{equation} \label{eq:domain} 
\sum _{i=r} ^{d}\Re(s_i)> d-r+1,\quad \quad i=1,2,\dots,d.
\end{equation}
Moreover, it is not difficult to
see that, for the domain of absolute convergence, it does not matter whether we
sum the summand on the right-hand side of \eqref{eq:multizeta} over
$0<n_1<\dots<n_d$ or over $0<n_1\le \dots\le n_d$. Therefore, the
domain described by the inequalities \eqref{eq:domain} is at the same
time the domain of absolute convergence of $\zeta(s_d,s_{d-1},\dots,s_1)$. 
That is, one has to add the conditions \eqref{eq:domain} for $i=2,\dots,d-1$ 
to Zhao's two conditions to obtain a complete description of the domain of
absolute convergence.
As a matter of fact, all the arguments given in the 
proof of Proposition~1 in \cite{ZhaoAB} are correct. However, 
it is only the case
$d=2$ which is carried out in detail (in which case there are no missing
conditions), and therefore the additional $d-2$ conditions are overlooked.
\end{Remark}

\begin{proof}[Proof of Proposition~\ref{prop:even}] (i) 
We consider first the left-hand side of
Andrews' identity~\eqref{eq:andrews}. We write the hypergeometric
series as a sum over $k$. Let $S_k$ denote the $k$-th summand.
Since for $N\ge k>\vert a\vert$ we have
$$
\left\vert\frac{(-N)_{k}}{(1+a+N)_{k}} 
\right\vert\le 
 \frac {(N-k+1)\cdots(N-1)N} 
{(N+1-\vert a\vert)(N+2-\vert a\vert)\cdots (N+k-\vert a\vert)}\le 1 ,
$$ 
and since for $k>N$ we have $(-N)_{k}=0$,
the modulus of $(-N)_{k}/(1+a+N)_{k}$ is bounded above by a constant
for {\it all\/} $k=0,1,\dots$.
Hence, using Lemma~\ref{lem:2}, we obtain that
$$
\vert S_k\vert\le 
D_{11}\cdot(k+1)^{-E-1},
$$
where $D_{11}$ is some constant independent of $k$, 
and where $E$ is the left-hand side of
\eqref{eq:cond1}. Since, by \eqref{eq:cond1}, we have $E>0$, the
absolutely convergent series $
\sum _{k=0} ^{\infty}D_{11}\cdot(k+1)^{-E-1}$ dominates the
hypergeometric series on the left-hand side of \eqref{eq:andrews}
term-wise. Thus, by Lebesgue's dominated convergence theorem, we
may perform its limit as $N\to\infty$ term-wise. This term-wise limit 
is exactly the left-hand side of \eqref{eq:Feven}.

Now we consider the right-hand side of \eqref{eq:andrews}. 
We need to temporarily assume that
\begin{equation} \label{eq:temp}
 \Re(a-b_{s+1}-c_{s+1})>0.
\end{equation}
(This is slightly stronger than \eqref{eq:cond2} with $r=s+1$.)
Writing $A$ for $a-b_{s+1}-c_{s+1}$,
for any non-negative integer $K\le N$ we have
$$
\left\vert\frac{(-N)_{K}}
{(b_{s+1}+c_{s+1}-a-N)_{K}}
\right\vert \le
\frac {N(N-1)\cdots(N-K+1)} {(N+\Re(A))
(N+\Re(A)-1)\cdots(N+\Re(A)-K+1)}
\le 1 ,
$$
and since for $K>N$ we have $(-N)_{K}=0$,
the modulus of $(-N)_{K}/(b_{s+1}+c_{s+1}-a-N)_{K}$ 
is bounded above by a constant
for {\it all\/} $K=0,1,\dots$. 
Thus, again using Lemma~\ref{lem:2}, the modulus of the summand indexed by
$k_1,k_2,\dots,k_s$ on the right-hand side of \eqref{eq:andrews} is
bounded above by
\begin{equation} \label{eq:expr} 
D_{12} 
\prod _{j=1} ^{s}(k_j+1)^{\Re(a-b_j-c_j)}\,
(k_1+\dots+k_j+1)^{\Re(b_j+c_j+b_{j+1}+c_{j+1}-2(a+1))},
\end{equation}
for some constant $D_{12}$ independent of the summation indices. Now we
apply Lemma~\ref{lem:4} with
$E_j=\Re(a-b_j-c_j)$ and $F_j=\Re(b_j+c_j+b_{j+1}+c_{j+1}-2(a+1))$.
This is indeed justified since, for this choice of parameters, the set
of conditions \eqref{eq:bed} is exactly the set \eqref{eq:cond2}. 
Hence, the sum of the expression \eqref{eq:expr} over all
$k_1,\dots,k_s\ge0$ converges.
Another application of Lebesgue's dominated convergence
theorem then implies that we may perform the limit of the multiple
sum on the right-hand side of \eqref{eq:andrews} as $N\to\infty$
term-wise. Together with the fact that
$$
\lim_{N\to+\infty}\frac{(1+a)_N\,(1+a-b_{s+1}-c_{s+1})_N}
{(1+a-b_{s+1})_N\,(1+a-c_{s+1})_N}
=\frac{\Ga(1+a-b_{s+1})\,\Ga(1+a-c_{s+1})}
{\Ga(1+a)\,\Ga(1+a-b_{s+1}-c_{s+1})}, 
$$ 
this establishes the identity \eqref{eq:Feven}, provided
\eqref{eq:temp} holds in addition to the conditions of the statement
of the proposition.

We can finally get rid of the restriction \eqref{eq:temp} by analytic
continuation. Indeed, by using arguments very similar to those above,
one can show that both sides of \eqref{eq:Feven} are analytic in
the parameters $a,b_1, \ldots, b_{s+1},c_1,
\ldots, c_{s+1}$ as long as \eqref{eq:cond1} and \eqref{eq:cond2} are
satisfied. In particular, in variation of Lemma~\ref{lem:2}, one would
use the fact that, for fixed complex numbers $\al$ and $\be$, 
there are constants $D_{13}$ and $D_{14}$ such that 
$$
 D_{13} \cdot(k+1)^{\Re(\al-\be)}
\log(k+2)\le\left\vert
\frac {\Ga(x+k)\,\psi(x+k)} {\Ga(\be+k)}\right\vert
\le D_{14} \cdot(k+1)^{\Re(\al-\be)}
\log(k+2)$$
for all non-negative integers $k$ and all complex numbers $x$ 
in a sufficiently small
neighbourhood of $\al$, say for $\vert
x-\al\vert<1$.
Here, $\psi(x)$ denotes the logarithmic derivative of $\Ga(x).$

(ii) In~\eqref{eq:Feven}, we first shift the parameters to $b_j\to
b_{j-1}$ and $c_{j}\to c_{j-1}$, and then we let $b_0\to+\infty$. The
same kind of argument as above then yields~\eqref{eq:Fodd}.
\end{proof}

\section{Multisum expansions of the Vasilyev-type integral $J_m$}\label{sec:int=andr}
The link between Andrews' identity and the Vasilyev-type integrals
$J_m$  
becomes apparent in the next proposition.
\begin{Proposition}\label{thm:zlobin}  
Let $z$, $a_0, a_1,\dots,a_{m}$, and $b_1,\dots,b_m$ be complex numbers
such that 
$\vert z\vert < 1$,  $\Re(a_0)>0$, $\Re(b_i)>\Re(a_i)>0$ 
for all $i=1,2,\dots,m$. 

{\em(i)} If $m=2s\ge 2$ is even, then
\begin{multline}\label{eq:Zlobin1}
J_{m}\left[\begin{matrix}
a_0,a_1,\dots,a_m\\b_1,\dots,b_m\end{matrix};z\right] =
\prod _{j=1} ^{2s}\frac {\Ga(a_{j})\,\Ga(b_{j}-a_{j})} {
\Ga(b_{j})}
\\
\times
\sum_{k_1, k_2,\ldots, k_{s} \ge 0} z^{k_1+\cdots+k_s}
\prod_{j=1}^s
\frac {(b_{2s-2j+2}-a_{2s-2j+2})_{k_j}} {k_j!}
\frac {(a_{2s-2j+1})_{k_1+\cdots +k_j}} 
{(b_{2s-2j+1})_{k_1+\cdots +k_j}}
\frac {(a_{2s-2j})_{k_1+\cdots +k_j}} 
{(b_{2s-2j+2})_{k_1+\cdots +k_j}}.
\end{multline}
This identity holds also for $z=1$ provided 
$\Re(b_1-a_1)\ge\Re(a_0)$, and provided \eqref{eq:bed} holds with
$E_j=\Re(b_{2s-2j+2}-a_{2s-2j+2}-1)$ and $F_j=\Re(a_{2s-2j}+a_{2s-2j+1}-
b_{2s-2j+1}-b_{2s-2j+2})$, $j=1,2,\dots,s$.

{\em(ii)} If $m=2s+1\ge 3$ is odd, then
\begin{multline}\label{eq:Zlobin2}
J_{m}\left[\begin{matrix}
a_0,a_1,\dots,a_m\\b_1,\dots,b_m\end{matrix};z\right]=
\prod _{j=1} ^{2s+1}\frac {\Ga(a_{j})\,\Ga(b_{j}-a_{j})} {
\Ga(b_{j})}
\cdot
\sum_{k_1,\ldots, k_{s+1} \ge 0} z^{k_1+\cdots+k_s}
\frac {(a_{2s+1})_{k_1}} {k_1!}
\frac {(a_{2s})_{k_1}} 
{(b_{2s+1})_{k_1}}\\
\times
\prod_{j=2}^{s+1}
\frac {(b_{2s-2j+4}-a_{2s-2j+4})_{k_j}} {k_j!}
\frac {(a_{2s-2j+3})_{k_1+\cdots +k_j}} 
{(b_{2s-2j+3})_{k_1+\cdots +k_j}}
\frac {(a_{2s-2j+2})_{k_1+\cdots +k_j}} 
{(b_{2s-2j+4})_{k_1+\cdots +k_j}}.\\
\end{multline}
This identity holds also for $z=1$ provided 
$\Re(b_1-a_1)\ge\Re(a_0)$, and provided~\eqref{eq:bed} holds with
$E_1=\Re(a_{2s+1}-1)$, $F_1=\Re(a_{2s}-b_{2s+1})$,
$E_j=\Re(b_{2s-2j+4}-a_{2s-2j+4}-1)$, and $F_j=\Re(a_{2s-2j+2}+a_{2s-2j+3}-
b_{2s-2j+3}-b_{2s-2j+4})$, $j=2,3,\dots,s+1$.
\end{Proposition}

\begin{proof} For $m\ge 2$, we denote by $Q_m(x_1, \ldots, x_{m} ; z)$ 
the nested expression in the denominator of the integrand 
in~\eqref{eq:Jdef}, that is
$$Q_m(x_1, \ldots, x_{m} ; z)=1-(1-(\cdots(1-x_m)x_{m-1})\cdots)x_1z.$$

(i) We prove the claim by induction on $m$. For $m=0$, the (empty)
integral $J_0\left[\begin{matrix} a_0\\ _{}
  \end{matrix};z\right]$ can be consistently interpreted as 1.
In order to do the induction step, we fix  $m=2s\ge 2$ 
and $z$ such that $\vert z\vert <1$. Then, trivially,    
\begin{multline*}
Q_{2s}(x_1, \ldots, x_{2s} ; z)=Q_{2s-2}(x_1, \ldots, x_{2s-2};z)-
zx_1 \cdots x_{2s-1}(1-x_{2s})\\=
Q_{2s-2}(x_1, \ldots ,x_{2s-2};z)\left(1-\frac{zx_1\cdots x_{2s-1}(1-x_{2s})}
{Q_{2s-2}(x_1, \ldots ,x_{2s-2};z)}\right),
\end{multline*}
where for $s=1$ the term $Q_0(-;z)$ has to be interpreted as 1.
Since for  $x_j\in[0,1]$, we have
$$
\left\vert\frac{zx_1\cdots x_{2s-1}(1-x_{2s})}
{Q_{2s-2}(x_1, \ldots ,x_{2s-2};z)} \right\vert \le \vert z\vert <1, 
$$
we may apply the binomial theorem to obtain
$$
\left(1-\frac{zx_1\cdots x_{2s-1}(1-x_{2s})}{Q_{2s-2}(x_1, \ldots
    ,x_{2s-2};z)}\right)^{-a_0}
=\sum_{k_1=0}^{\infty}z^{k_1}\,\frac{\Ga(a_0+k_1)}{\Ga(a_0)\Ga(k_1+1)}
\left(\frac{x_1\cdots x_{2s-1}(1-x_{2s})}
{Q_{2s-2}(x_1, \ldots ,x_{2s-2};z)}\right)^{k_1}.
$$
Hence
\begin{multline*}
J_{2s}\left[\begin{matrix}
a_0,a_1,\dots,a_{2s}\\b_1,\dots,b_{2s}\end{matrix};z\right]=
\int_{[0,1]^{2s}}\sum_{k_1=0}^{\infty}z^{k_1}\,
\frac{\Ga(a_0+k_1)}{\Ga(a_0)\Ga(k_1+1)} 
 x_{2s-1}^{a_{2s-1}+k_1-1}(1-x_{2s-1})^{b_{2s-1}-a_{2s-1}-1}\\
\cdot  x_{2s}^{a_{2s}-1}(1-x_{2s})^{b_{2s}-a_{2s}+k_1-1}
\frac{\prod_{j=1}^{2s-2} x_j^{a_j+k_1-1}(1-x_j)^{b_j-a_j-1}}
{Q_{2s-2}(x_1, \ldots ,x_{2s-2};z)^{k_1+a_0}}\, \dd x_1\cdots \dd x_{2s}.
\end{multline*}
The conditions on the parameters ensure that the integral
\begin{multline*}
\int_{[0,1]^{2s}}\sum_{k_1=0}^{\infty}\bigg\vert z^{k_1}\,
\frac{\Ga(a_0+k_1)}{\Ga(a_0)\Ga(k_1+1)} 
 x_{2s-1}^{a_{2s-1}+k_1-1}(1-x_{2s-1})^{b_{2s-1}-a_{2s-1}-1}\\
\cdot  x_{2s}^{a_{2s}-1}(1-x_{2s})^{b_{2s}-a_{2s}+k_1-1}
\frac{\prod_{j=1}^{2s-2} x_j^{a_j+k_1-1}(1-x_j)^{b_j-a_j-1}}
{Q_{2s-2}(x_1, \ldots ,x_{2s-2};z)^{a_0+k_1}}\bigg\vert\, \dd x_1\cdots \dd x_{2s}
\end{multline*}
is convergent. Thus, we can  
exchange the integral and
the summation, and, using the beta integral evaluation and some 
straightforward simplifications, we obtain
\begin{multline}
J_{2s}\left[\begin{matrix}
a_0,a_1,\dots,a_{2s}\\b_1,\dots,b_{2s}\end{matrix};z\right]
= \frac{\Ga(a_{2s})\,\Ga(a_{2s-1})\,\Ga(b_{2s}-a_{2s})\,
\Ga(b_{2s-1}-a_{2s-1})}
 {\Ga(b_{2s})\,\Ga(b_{2s-1})}\\
\cdot\sum_{k_1=0}^{\infty}z^{k_1}\,
\frac{(b_{2s}-a_{2s})_{k_1}\,(a_{2s-1})_{k_1}\,(a_0)_{k_1}}
{k_1!\,(b_{2s})_{k_1}\,(b_{2s-1})_{k_1}} J_{2s-2}\left[\begin{matrix}
a_0+k_1,a_1+k_1,\dots,a_{2s-2}+k_1\\b_1+k_1,\dots,b_{2s-2}+k_1\end{matrix};z\right].
\label{eq:zlobeven}
\end{multline}
If we substitute the induction hypothesis for $J_{2s-2}$, we arrive 
exactly at \eqref{eq:Zlobin1}.

We now perform the limit $z\to1$: 
since the conditions on the parameters  guarantee 
that the integral $J_{2s}$ is absolutely convergent for $z=1$,  
dominated convergence implies that one can
interchange limit and integral. Similarly, if we put $z=1$ in the above
multiple sum, then the conditions on the
parameters allow us to apply Lemmas~\ref{lem:2} and \ref{lem:4} and to
conclude that it converges absolutely. Thus, again, dominated
convergence implies that we may interchange limit and summation.
As a result, Case (i) of Proposition~\ref{thm:zlobin} is now completely proved.

\medskip

(ii) 
We do not provide all the details for the case where $m$ is odd, 
$m=2s+1\ge 3$,
since this case can be treated in a rather similar manner as
the case where $m$ is even. 
A main difference, however, is that, to get started,
we use the alternative identity  
$$
Q_{2s+1}(x_1, \ldots, x_{2s+1} ; z)=
Q_{2s}(x_1, \ldots,x_{2s};z)\left(1-\frac{zx_1 \cdots x_{2s+1}}
{Q_{2s}(x_1, \ldots,x_{2s};z)}\right),
$$
which implies the expansion
\begin{multline}
J_{2s+1}\left[\begin{matrix}
a_0,a_1,\dots,a_{2s+1}\\b_1,\dots,b_{2s+1}\end{matrix};z\right]
\\= \frac{\Ga(a_{2s+1})\,\Ga(b_{2s+1}-a_{2s+1})}
 {\Ga(b_{2s+1})}\cdot\sum_{k=0}^{\infty}z^k  \,
\frac{(a_0)_k\,(a_{2s+1})_k}
{k!\,(b_{2s+1})_k}
J_{2s}\left[\begin{matrix}
a_0+k,a_1+k,\dots,a_{2s}+k\\b_1+k,\dots,b_{2s}+k\end{matrix};z\right].
\label{eq:zlobodd}
\end{multline}
At this point, we substitute the multiple
series~\eqref{eq:Zlobin1} 
for $J_{2s}$, and after some simple manipulations
we arrive at \eqref{eq:Zlobin2}. 
\end{proof}

\begin{Remark}
Both of the recursive formulas \eqref{eq:zlobeven} and \eqref{eq:zlobodd}
appear already earlier in the article \cite{zlo} of Zlobin which we
mentioned in the Introduction.
He used them to express the
integrals $J_{m}$ in terms of another family of integrals, like those
considered by Sorokin in~\cite{so, so1}. 
\end{Remark}

\section{Proof of Theorem~\ref{thm:zudgene}}\label{sec:proofvaszud}
We are now in the position to prove Zudilin's theorem, by putting together
the identities in Proposition~\ref{thm:zlobin} and~\ref{prop:even}.
Because of the use of Proposition~\ref{thm:zlobin} when $z=1$, 
we shall need to temporarily
impose stronger conditions on the parameters than required by the
assertion of the theorem. We shall do this without mention. 
One gets rid of these restrictions at the
end by analytic continuation.

Let first $m$ be even, $m=2s$.
We apply Proposition~\ref{thm:zlobin}, Eq.~\eqref{eq:Zlobin1}, with 
$a_{j-1}=h_j$, $j=1,2,\dots,2s+1$,
$b_{j}=1+h_0-h_{j+2}$, $j=1,2,\dots,2s$. Thus, using
\eqref{eq:Zlobin1}, we express the integral on the left-hand side of
\eqref{eq:Zudilin} in terms of a multiple sum. If we subsequently
apply the identity \eqref{eq:Feven} with
$b_j=h_{2s-2j+4}$, $c_j=h_{2s-2j+3}$ 
for  $j=1,2, \dots, s+1$,
to the multiple sum, then we
arrive at the very-well-poised hypergeometric series on the right-hand
side of \eqref{eq:Zudilin}.

Similarly, if $m$ is odd, $m=2s+1$, then we 
apply Proposition~\ref{thm:zlobin}, Eq.~\eqref{eq:Zlobin2}, with 
$a_{j-1}=h_j$, $j=1,2,\dots,2s+2$,
$b_{j}=1+h_0-h_{j+2}$, $j=1,2,\dots,2s+1$. Thus, using
\eqref{eq:Zlobin2}, we express the integral on the left-hand side of
\eqref{eq:Zudilin} in terms of a multiple sum. If we subsequently
apply the identity \eqref{eq:Fodd} with $s$ replaced by $s+1$, 
$b_j=h_{2s-2j+4}$, $j=1,2, \dots, s+1$, $c_j=h_{2s-2j+3}$ 
for  $j=0,1, \dots, s+1$,
to the multiple sum, then we
arrive at the very-well-poised hypergeometric series on the right-hand
side of \eqref{eq:Zudilin}.


\section*{Acknowledgement}

We thank Wadim Zudilin for an attentive reading of an earlier version of the
paper, and, in particular, for several useful 
suggestions for improvement of the exposition.

\def\refname{Bibliography}

\end{document}